\documentclass[12pt]{article}

\usepackage{amsmath, epsfig, cite, setspace}
\usepackage{amssymb}
\usepackage{amsfonts}
\usepackage{latexsym}
\usepackage{amsthm}

\makeatletter
\renewcommand{\@seccntformat}[1]{{\csname the#1\endcsname}.\hspace{.5em}}
\makeatother

\newtheorem{thm}{Theorem}[section]
\newtheorem{prop}[thm]{Proposition}

\newtheorem{conj}[thm]{Conjecture}

\renewcommand{\qed}{\hfill$\Box$\medskip}
\renewcommand{\thefootnote}{*}

\numberwithin{equation}{section}

\begin{document}

\begin{center}
{\large\bf Congruence properties for the trinomial coefficients}
\end{center}

\vskip 2mm \centerline{Moa Apagodu$^1$ and Ji-Cai Liu$^2$}
\begin{center}
{\footnotesize $^1$Department of Mathematics, Virginia Commonwealth University, Richmond, VA 23284\\
{\tt mapagodu@vcu.edu} \\[10pt]

$^2$Department of Mathematics, Wenzhou University, Wenzhou 325035, People's Republic of China\\
{\tt jcliu2016@gmail.com  } }
\end{center}


\vskip 0.7cm \noindent{\bf Abstract.}
In this paper, we state and prove some congruence properties for the trinomial coeficients,
one of which is similar to the Wolstenholme's theorem.



\renewcommand{\thefootnote}{**}

\section{Introduction}
In 1819, Babbage \cite{babbage-epj-1819} showed for any odd prime $p$,
\begin{align}
{2p-1\choose p-1}\equiv 1\pmod{p^2}.\label{bab}
\end{align}
In 1862, Wolstenholme \cite{wolstenholme-qjpam-1862} proved that the above congruence holds modulo $p^3$ for any prime $p\ge 5$,
which is known as the famous Wolstenholme's theorem. It is well-known that Wolstenholme's theorem is a fundamental congruence in combinatorial number theory. We refer to \cite{mestrovic-a-2011} for various extensions of Wolstenholme's theorem.

In the past few years, ($q$-)congruences for sums of binomial coefficients have attracted the attention of many researchers (see, for instance, \cite{apagodu-ijnt-2018,az-amm-2017,guo-ijnt-2018,gs-2018,gz-aam-2010,gz-am-2019,liu-jdea-2016,liu-rm-2018,st-aam-2010,st-ijnt-2011}).
In 2011, Sun and Tauraso \cite{st-ijnt-2011} proved that for any prime $p\ge 5$,
\begin{align}
&\sum_{k=0}^{p-1}{2k\choose k}\equiv \left(\frac{p}{3}\right)\pmod{p^2},\label{st-1}\\
&\sum_{k=0}^{p-1}{2k\choose k}\frac{1}{k+1}\equiv\frac{3}{2}\left(\frac{p}{3}\right)- \frac{1}{2}\pmod{p^2},\label{st-2}
\end{align}
where $\left(\frac{\cdot}{p}\right)$ denotes the Legendre symbol.
Note that ${2k\choose k}\frac{1}{k+1}$ is the $n$th Catalan number $C_n$,
which play an important role in various counting problems. Extensions of \eqref{st-1} and \eqref{st-2} have been established in \cite{az-amm-2017,liu-jdea-2016}.

In 2018, the first author \cite{apagodu-ijnt-2018} conjectured two congruences on sums of the super Catalan numbers (named by Gessel \cite{gessel-jsc-1992}):
\begin{align*}
\sum_{i=0}^{p-1}\sum_{j=0}^{p-1}\frac{{2i\choose i}{2j\choose j}}{{i+j\choose i}}\equiv \left(\frac{p}{3}\right)\pmod{p},\\[5pt]
\sum_{i=0}^{p-1}\sum_{j=0}^{p-1}(3i+3j+1)\frac{{2i\choose i}{2j\choose j}}{{i+j\choose i}}\equiv -7\left(\frac{p}{3}\right)\pmod{p},
\end{align*}
which were confirmed by the second author \cite{liu-rm-2018}.

In this paper, we will study the congruence properties for the trinomial coefficients.
Here we consider the coefficients of the trinomial:
$$
(1+x+x^{-1})^n=\sum_{j=-n}^{n}\left( n \choose j \right)x^j.
$$
Two immediate consequences of this definition are
$$
\left( n \choose j \right) =\left(n \choose -j \right),
$$
and
$$
\left( n \choose j \right) =\left(n-1 \choose j-1 \right)+\left(n-1 \choose j \right)+\left(n-1 \choose j+1 \right).
$$
It is not hard to show that (see \cite{sills-b-2018}) the identities
\begin{align}
\left(n \choose j \right)=\sum_{k=0}^{n}
{n \choose k} { n-k \choose k+j},\label{new-1}
\end{align}
and
$$
\left(n \choose j \right)=\sum_{k=0}^{n}
(-1)^k{n \choose k} { 2n-2k \choose n-k-j}.
$$

We first prove a congruence for the trinomial coefficients, which is similar to the
Wolstenholme's theorem.
\begin{thm}\label{t-1}
For any prime $p\ge 5$, we have
\begin{align}
\left(2p\choose p\right)\equiv 2\pmod{p^2}.\label{new-2}
\end{align}
\end{thm}

The second result consists of the following two congruences on single sums of trinomial coefficients.
\begin{thm}\label{t-2}
For any prime $p\ge 5$, we have
\begin{align}
&\sum_{j=0}^{p}\left(p\choose j\right)\equiv \frac{1}{2}\left(1+3^p\right)\pmod{p^2},\label{a-1}\\[5pt]
&\sum_{j=0}^{p-1}\left(p-1\choose j\right)\equiv \frac{1}{2}\left(1+\left(\frac{p}{3}\right)\right)\pmod{p}.\label{a-2}
\end{align}
\end{thm}

The third aim of the paper is to establish a congruence on double sums of trinomial coefficients.
\begin{thm}\label{t-3}
For any prime $p\ge 5$ and integer $j$ with $0<j<p$, we have
\begin{align}
&\sum_{k=0}^{p-1}\left(k\choose j\right)\equiv \frac{(-1)^j+1}{2}\cdot (-1)^{\frac{p-j-1}{2}}
\pmod{p}\label{a-3},\\[5pt]
&\sum_{m=0}^{p-1}\sum_{n=0}^{p-1}\left({m\choose n}\right)\equiv \frac{1}{2}\left((-1)^{\frac{p-1}{2}}+1\right)\pmod{p}.\label{a-4}
\end{align}
\end{thm}

The rest of this paper is organized as follows.
We shall prove Theorems \ref{t-1}, \ref{t-2} and \ref{t-3} in Sections 2, 3 and 4, respectively.
An open problem on $q$-congruence is proposed in the last section for further research.

\section{Proof of Theorem \ref{t-1}}
By \eqref{new-1}, we have
\begin{align*}
\left(2p\choose p\right)=\sum_{k=0}^{\frac{p-1}{2}}
{2p \choose k} { 2p-k \choose k+p}.
\end{align*}
For $1\le k\le \frac{p-1}{2}$, we have
\begin{align*}
{2p\choose k}=\frac{2p(2p-1)\cdots(2p-k+1)}{k!}\equiv \frac{2(-1)^{k-1}p}{k}\pmod{p^2},
\end{align*}
and so
\begin{align}
\left(2p\choose p\right)\equiv {2p\choose p}+2p\sum_{k=1}^{\frac{p-1}{2}}
\frac{(-1)^{k-1}}{k} {2p-k \choose k+p}\pmod{p^2}.\label{new-5}
\end{align}

Furthermore,
\begin{align*}
{2p-k \choose k+p}
&=\frac{(2p-k)(2p-k-1)\cdots(p+1)\cdot (p-1)(p-2)\cdots (p-2k+1)}{(p+1)(p+2)\cdots(p+k)\cdot(p-1)!}\\[5pt]
&\equiv -\frac{(p-k)!(2k-1)!}{k!(p-1)!}\pmod{p}\\[5pt]
&=-\frac{(2k-1)!}{k!(p-1)(p-2)\cdots(p-k+1)}\\[5pt]
&\equiv \frac{(-1)^k(2k-1)!}{k!(k-1)!}\pmod{p}\\[5pt]
&=\frac{(-1)^k}{2}{2k\choose k}.
\end{align*}
It follows from the above and \eqref{new-5} that
\begin{align}
\left(2p\choose p\right)\equiv 2{2p-1\choose p-1}-p\sum_{k=1}^{\frac{p-1}{2}}
{2k\choose k}\frac{1}{k}\pmod{p^2}.\label{new-3}
\end{align}
By \cite[Theorem 1.3]{st-aam-2010}, we have
\begin{align}
\sum_{k=1}^{\frac{p-1}{2}}\frac{{2k\choose k}}{k}\equiv 0\pmod{p}.\label{new-4}
\end{align}
Combining \eqref{bab}, \eqref{new-3} and \eqref{new-4}, we complete the proof of \eqref{new-2}.

\section{Proof of Theorem \ref{t-2}}
{\noindent\it Proof of \eqref{a-1}.}
We begin with the following identity (see {\tt http://oeis.org/A027914})
\begin{align}
\sum_{j=0}^n\left(n\choose j\right)=\frac{1}{2}\left(3^n+\sum_{k=0}^n{n\choose k}{n-k\choose k}\right).\label{b-1}
\end{align}
Letting $n=p$ in the above gives
\begin{align*}
\sum_{j=0}^{p}\left(p\choose j\right)&=\frac{1}{2}\left(3^{p}+\sum_{k=0}^{\frac{p-1}{2}}{p\choose k}{p-k\choose k}\right).
\end{align*}
Note that for $1\le k\le \frac{p-1}{2}$,
\begin{align*}
{p\choose k}{p-k\choose k}&=\frac{p(p-1)\cdots(p-k+1)(p-k)(p-k-1)\cdots(p-2k+1)}{k!^2}\\
&\equiv -\frac{p}{2k}{2k\choose k}\pmod{p^2}.
\end{align*}
Thus,
\begin{align*}
\sum_{j=0}^{p}\left(p\choose j\right)\equiv \frac{1}{2}\left(3^{p}+1-\frac{p}{2}\sum_{k=1}^{\frac{p-1}{2}}\frac{{2k\choose k}}{k}\right)\pmod{p^2}.
\end{align*}
It follows from the above and \eqref{new-4} that
\begin{align*}
\sum_{j=0}^{p}\left(p\choose j\right)\equiv \frac{1}{2}\left(3^{p}+1\right)\pmod{p^2},
\end{align*}
as desired.
\qed

{\noindent \it Proof of \eqref{a-2}.}
Letting $n=p-1$ in \eqref{b-1}, we obtain
\begin{align*}
\sum_{j=0}^{p-1}\left(p-1\choose j\right)&=\frac{1}{2}\left(3^{p-1}+\sum_{k=0}^{\frac{p-1}{2}}{p-1\choose k}{p-1-k\choose k}\right).
\end{align*}
Note that for $0\le k\le \frac{p-1}{2}$,
\begin{align*}
{p-1\choose k}\equiv (-1)^k \pmod{p},
\end{align*}
and
\begin{align*}
{p-1-k\choose k}=\frac{(p-1-k)(p-2-k)\cdots(p-2k)}{k!}\equiv (-1)^k{2k\choose k}\pmod{p}.
\end{align*}
By the above two congruences and Fermat's little theorem, we have

\begin{align}
\sum_{j=0}^{p-1}\left(p-1\choose j\right)\equiv \frac{1}{2}\left(1+\sum_{k=0}^{p-1}{2k\choose k}\right)\pmod{p}.\label{b-2}
\end{align}
Then the proof of \eqref{a-2} follows from \eqref{st-1} and \eqref{b-2}.
\qed

\section{Proof of Theorem \ref{t-3}}
{\noindent \it Proof of \eqref{a-3}.}
Exchanging the summation order, we get
\begin{align*}
\sum_{k=0}^{p-1}\left(k\choose j\right)&=\sum_{k=0}^{p-1}\sum_{i=0}^k{k\choose i}{k-i\choose i+j}\\
&=\sum_{i=0}^{p-1}\sum_{k=i}^{p-1}{k\choose i}{k-i\choose i+j}.
\end{align*}
Since
\begin{align*}
{k\choose i}{k-i\choose i+j}={2i+j\choose i}{k\choose 2i+j},
\end{align*}
we have
\begin{align*}
\sum_{k=0}^{p-1}\left(k\choose j\right)
&=\sum_{i=0}^{p-1}{2i+j\choose i}\sum_{k=i}^{p-1}{k\choose 2i+j}\\
&=\sum_{i=0}^{p-1}{2i+j\choose i}{p\choose 2i+j+1},
\end{align*}
where we have utilized the identity (proved by induction):
\begin{align*}
\sum_{k=0}^n{k\choose m}={n+1\choose m+1}.
\end{align*}
Note that for $1\le k\le p-1$,
\begin{align*}
{p\choose k}\equiv 0\pmod{p}.
\end{align*}

\noindent If $j$ is odd, then
\begin{align*}
\sum_{k=0}^{p-1}\left(k\choose j\right)\equiv 0\pmod{p}.
\end{align*}

\noindent If is $j$ is even, then
\begin{align*}
\sum_{k=0}^{p-1}\left(k\choose j\right)\equiv {p-1\choose \frac{p-j-1}{2}}\equiv (-1)^{\frac{p-j-1}{2}}\pmod{p}.
\end{align*}
This completes the proof of \eqref{a-3}.
\qed

{\noindent \it Proof of \eqref{a-4}.}
By \eqref{a-3}, we have
\begin{align*}
\sum_{m=0}^{p-1}\sum_{n=0}^{p-1}\left({m\choose n}\right)&\equiv
\sum_{n=0}^{p-1}\frac{(-1)^n+1}{2}\cdot (-1)^{\frac{p-n-1}{2}}\pmod{p}\\
&=\sum_{n=0}^{\frac{p-1}{2}}(-1)^{\frac{p-2n-1}{2}}\\
&=\frac{1}{2}\left((-1)^{\frac{p-1}{2}}+1\right),
\end{align*}
as claimed.
\qed

\noindent {\bf Remark.} Theorems \ref{t-2} and \ref{t-3} can also be established by using the Method of the first author and Zeilberger \cite{az-amm-2017}.

\section{Concluding remarks}
We have three $q$-analogs corresponding to the trinomial coefficients as given in \cite{sills-b-2018}, namely,
$$
T_1(n,j,q):=\sum_{k=0}^{n}
q^{k(k+j)}{n \choose k}_q { n-k \choose k+j}_q,
$$

$$
T_2(n,j,q):=\sum_{k=0}^{n}
(-1)^k{n \choose k}_{q^2} { 2n-2k \choose n-k-j}_q,
$$

$$
T_3(n,j,q):=\sum_{k=0}^{n}
(-q)^k{n \choose k}_{q^2} { 2n-2k \choose n-k-j}_q,
$$
where the $q$-binomial coefficients ${n\choose k}_q$ are defined as
$$
{n\choose k}_q
=\begin{cases}
\displaystyle\frac{(1-q^n)(1-q^{n-1})\cdots (1-q^{n-k+1})}{(1-q)(1-q^2)\cdots(1-q^k)}, &\text{if $0\leqslant k\leqslant n$},\\[10pt]
0,&\text{otherwise.}
\end{cases}
$$

It is not hard to prove the following $q$-congruences.
\begin{prop}\label{p-1}
For any odd prime $p$ and integer $1\le s\le 3$, we have
\begin{align*}
\sum_{j=0}^{p-1}T_s(p,j,q)\equiv 1 \pmod{[p]_q},
\end{align*}
where the $q$-integers are given by $[n]_q=(1-q^n)/(1-q)$.
\end{prop}
\noindent The proof of Proposition \ref{p-1} is trivial and left to the interested reader.

In 1999, Andrews \cite{andrews-dm-1999} established an interesting $q$-analog of Babbage's congruence \eqref{bab}:
\begin{align*}
{2p-1\choose p-1}_q\equiv q^{\frac{p(p-1)}{2}}\pmod{[p]_q^2},
\end{align*}
for any odd prime $p$. It is natural to ask whether the congruence \eqref{new-2} possesses a
$q$-analog. For convenience sake, let
$$
\left({n\choose j}\right)_q=T_1(n,j,q).
$$
Numerical calculation suggests the following $q$-congruence, and we
propose this conjecture for further research.
\begin{conj}
For any prime $p\ge 5$, we have
\begin{align}
\left(2p\choose p\right)_q \equiv \left(2\left\lfloor \frac{p+3}{6}\right\rfloor+p\right)(q^p-1)+2 \pmod{[p]_q^2},\label{conj}
\end{align}
where $\lfloor x \rfloor$ denotes the integral part of real $x$.
\end{conj}
\noindent Letting $q\to 1$ in \eqref{conj}, we are led to \eqref{new-2}.

\vskip 3mm \noindent{\bf Acknowledgments.}
The first author would like to thank Professor Dennis Stanton of the School of Mathematics, University of Minnesota for indoctrinating him into the $q$-series during his visit in the fall 2019. The second author was supported by the National Natural Science Foundation of China (grant 11801417).

\end{document}